\documentclass[11pt,a4paper,leqno]{article}
\usepackage{a4wide}
\usepackage{latexsym}
\usepackage{amsmath}
\usepackage{amssymb}
\usepackage{stackrel}  
\usepackage{color} 
\usepackage{graphicx}
\usepackage[linesnumbered,ruled,vlined]{algorithm2e}

\usepackage{authblk}

\usepackage{hyperref}

\newtheorem{defin}{Definition}
\newtheorem{lemma}{Lemma}
\newtheorem{prop}{Proposition}
\newtheorem{theo}{Theorem}
\newtheorem{corol}{Corollary}

\newenvironment{proof}{\medskip\par\noindent{\bf Proof}}{\hfill $\Box$
\medskip\par}

\newcommand{\C}{\mathbb{C}}

\newcommand{\R}{\mathbb{R}}
\newcommand{\Z}{\mathbb{Z}}

\begin{document}
\title{On parametric $0$-Gevrey asymptotic expansions in two levels for some linear partial $q$-difference-differential equations}

\author[1]{Alberto Lastra}
\author[2]{St\'ephane Malek}
\affil[1]{Universidad de Alcal\'a, Dpto. F\'isica y Matem\'aticas, Alcal\'a de Henares, Madrid, Spain. {\tt alberto.lastra@uah.es}}
\affil[2]{University of Lille, Laboratoire Paul Painlev\'e, Villeneuve d'Ascq cedex, France. {\tt stephane.malek@univ-lille.fr}}


\date{}

\maketitle
\thispagestyle{empty}
{ \small \begin{center}
{\bf Abstract}
\end{center}

A novel asymptotic representation of the analytic solutions to a family of singularly perturbed $q-$difference-differential equations in the complex domain is obtained. Such asymptotic relation shows two different levels associated to the vanishing rate of the domains of the coefficients in the formal asymptotic expansion. On the way, a novel version of a multilevel sequential Ramis-Sibuya type theorem is achieved.


\smallskip

\noindent Key words: q-Gevrey asymptotic expansions; singularly perturbed; formal power series. 2020 MSC: 35C10, 35R10, 35C15, 35C20.
}
\bigskip \bigskip

\section{Introduction}

This work is devoted to construct asymptotic formal expansions related to sectorial analytic ones to certain family of $q-$difference$-$differential equations in the complex domain. 

Let $q>1$ and let $\sigma_{q,t}$ be the dilation operator on $t$ variable defined by $\sigma_{q,t}
f(t)=f(qt)$. The elements of the family of functional equations under study are of the following form: 
\begin{equation}\label{epralintro}
Q(\partial_z)\sigma_{q,t}u(t,z,\epsilon)=F(t,z,\epsilon,\sigma_{q,t}.\partial_z)u(t,z,\epsilon)+\sigma_{q,t}f(t,z,\epsilon),
\end{equation}
where $Q(X)$ is a polynomial, $F(x_1,\ldots,x_5)$ (and also $f(x_1,x_2,x_3)$) is a holomorphic function with respect to $(x_1,x_2,x_3)$ on the product of a sector with vertex at the origin times a horizontal strip  times a neigborhood of the origin; algebraic with respect to $x_4$, and of polynomial nature in its last variable. The element $\epsilon$ acts in the equation as a small complex perturbation parameter. Let us denote $\mathcal{T}\times H\times D$ the domain of definition of $F$ with respect to its three first variables. The precise conditions on the elements involved in the definition of equation (\ref{epralintro}) are determined in Section~\ref{sec1}.

It is worth mentioning that we have distinguished two terms in equation (\ref{epralintro}), which turn out to be the leading terms, and are of the form
$$(\epsilon t)^{d_{D_j}}\sigma_{q,t}^{\frac{d_{D_j}}{k_j}+1}R_{D_j}(\partial_z)u(t,z,\epsilon),\qquad j=1,2,$$
for some integers $d_{D_j}\ge1$ and polynomials $R_{D_j}$, where $k_j>0$ are given real numbers. In other words, (\ref{epralintro}) reads as follows:
\begin{multline*}
 Q(\partial_z)\sigma_{q,t}u(t,z,\epsilon)=\sum_{j=1,2}(\epsilon t)^{d_{D_j}}\sigma_{q;t}^{\frac{d_{D_j}}{k_j}+1}R_{D_j}(\partial_z)u(t,z,\epsilon)\\
+G(t,z,\epsilon,\sigma_{q,t},\partial_z)u(t,z,\epsilon)+\sigma_{q,t}f(t,z,\epsilon),
\end{multline*}
for $G(x_1,\ldots,x_5)$ with the same nature as that of $F$, and polynomial coefficients in $x_4$. The first two terms, and more precisely the numbers $k_1,k_2$ determine the asymptotic behavior of the formal solution when the perturbation parameter $\epsilon$ approaches the origin. 

In the previous work~\cite{drelasmal}, the authors have constructed a family of analytic solutions to (\ref{epralintro}), say $(u_p(t,z,\epsilon))_{0\le p\le \varsigma-1}$ for some $\varsigma\ge2$, which are defined in sectorial regions with respect to the perturbation parameter (see Theorem~\ref{teo1}). Such sectorial regions cover a punctured disc at the origin, so that the difference of two consecutive solutions (consecutive in the sense that they are solutions defined on consecutive sectors in $\epsilon$) can be estimated in the intersection of these domains in order to apply a Ramis-Sibuya type theorem involving different levels. Indeed, the main result in~\cite{drelasmal} is a consequence of that previous result, and states a decomposition of the analytic  solution and its formal expansion to the main problem as the sum of two terms. In that decomposition one of the terms related to the analytic solutions admits one of the elements of the splitting of the formal expansion as its common $q-$Gevrey asymptotic expansion of order $k_1$ on the corresponding sector in $\epsilon$ whereas the other part of the splitting in the analytic solution admits the remaining part of the formal expansion as its $q-$Gevrey asymptotic expansion of order $k_2$ in that same sector. The details of this result are recalled in Theorem~\ref{teopral0}.

This work is a continuation of~\cite{drelasmal}. We provide a different asymptotic representation of the analytic solutions to (\ref{epralintro}) following the novel concept of asymptotic expansion introduced by H. Tahara in~\cite{tahara24} (see Definition~\ref{defi319}). Such representation is motivated by the shape of the formal solutions to $q$-difference equations, observed in the recent work~\cite{carrillolastra}. This asymptotic scheme is based on the approximation of an analytic function in $(t,\epsilon)$ by a power series of the form
\begin{equation}\label{e109}
\sum_{p\ge0}\varphi_p(t)\epsilon^p,
\end{equation}
where $(\varphi_p)_{p\ge0}$ is a family of holomorphic functions, with $\varphi_p$ being a function which is defined on some open domain $U_p\subseteq\C$ in such a way that $U_{p+1}\subseteq U_{p}$ for all $p\ge0$, and $\hbox{diam}(U_p)\to 0$ when $p\to\infty$. Here, $\hbox{diam}(\cdot)$ stands for the diameter of a set, i.e. the supreme of the distances of two of its points. The power series (\ref{e109}) determines a formal power series in $\epsilon$ with null radius of convergence. The concept of Gevrey and $q-$Gevrey asymptotic expansions and orders can be adapted from the classical settings to such approximation scheme. In this sense, in this work we deal with $0-$Gevrey asymptotic expansions relative to the sequence $(\hbox{diam}(U_p))_{p\to\infty}$ of a function, i.e. the upper estimates of the approximation of the analytic function and the truncation of the formal power series being of geometric nature.

The main result of the present work states the existence of a splitting of the analytic and the formal power series expansions in the form
$$u_p(t,z,\epsilon)=\mathfrak{u}_{p}^{1}(t,z,\epsilon)+\mathfrak{u}_{p}^{2}(t,z,\epsilon).$$ 
Each of the terms in the splitting  admits a corresponding element in the splitting of the formal power series expansion 
$$\hat{u}(t,z,\epsilon)=\hat{\mathfrak{u}}_1(t,z,\epsilon)+\hat{\mathfrak{u}}_2(t,z,\epsilon),$$
as its $0-$Gevrey asymptotic expansion relative to two different sequences of positive numbers, say $r_1=(r_{1,p})_{p\ge0}$ and $r_2=(r_{2,p})_{p\ge0}$, whose rate of decay to 0 differs one to each other (see Theorem~\ref{teopral}). For this purpose, we develop a multilevel sequential Ramis-Sibuya (RS) Theorem, which is the sequential counterpart of the multilevel version of Ramis-Sibuya Theorem which has already been applied in the literature. As a matter of fact, the splitting of the analytic and their formal expansions resembles that in the framework of multisummability results, as described in Theorem 50,~\cite{ba2} or in Section 7.5,~\cite{loday} with more details. 

The study of the asymptotic behavior of solutions to $q$-difference equations in the complex domain in the vicinity of an irregular singularity has been a topic of interest in the last decades. We mention the works~\cite{dre,lamasa0,ma11,tahara2,taya2,vizh,ya1}, among many others. There is also an increasing interest on the knowledge of sequences preserving summability and $q-$Gevrey asymptotics, which simplifies computations, such as the works~\cite{icmi,lami}. The applications of this theory are also of great importance and interest,~\cite{prrasp,prrasp1}. 

We have decided to include the main results in~\cite{drelasmal} for the sake of completeness, as some of the results in the present work heavily lean on the geometric construction of the solutions in that work, or the hypotheses made. More precisely, the paper is organized as follows. Section~\ref{sec1} is devoted to recall the main equation under study, the assumptions made on its elements, and the form of the analytic solutions which are constructed in the form of an inverse and $q-$Laplace transform. Some auxiliary results as Proposition~\ref{prop227} are needed in the sequel. The asymptotic representation of the analytic solutions described in~\cite{drelasmal} is recalled in Section~\ref{sec31}. In Section~\ref{sec32}, we obtain novel upper bounds for the difference of two consecutive analytic solutions in different situations (Proposition~\ref{prop293} and Proposition~\ref{prop269}). The work concludes with Section~\ref{secpral}, where a multilevel sequential Ramis-Sibuya (RS) Theorem is attained, and applied to achieve the main result of the work, Theorem~\ref{teopral}. The work concludes with an Appendix where we recall some facts on classical operators and functions involved in the achievement of the results in the paper.

\vspace{0.3cm}

\textbf{Notation:}

For every $q>1$, we denote the operator $\sigma_{q,t}$ defined by  $\sigma_{q,t}(f)(t)=f(qt)$, whenever it makes sense, when applied on $f(t)$. Given any $\delta>0$, we extend the previous definition to $\sigma_{q,t}^{\delta}(f)(t)=f(q^{\delta}t)$.

For $z_0\in\C$ and $r>0$, we write $D(z_0,r)$ for the open disc $\{z\in\C:|z-z_0|<r\}$, and $\overline{D}(z_0,r)$ for its closure. Given any sector of the complex domain $S$ with vertex at the origin, we say that $T$ is a proper subsector of $S$ and we write $T\prec S$ if $T$ is a bounded open sector with vertex at the origin, and $\overline{T}\setminus\{0\}\subseteq S$.

Let $\mathbb{E}$ be a Banach space of functions. Given an open set $U\subseteq\R$, $\mathcal{C}(U,\mathbb{E})$ stands for the set of continuous functions defined on $U$ with values in $\mathbb{E}$. We write $\mathbb{E}\{\epsilon\}$ for the set of holomorphic functions defined on some neighborhood of the origin (in the variable $\epsilon$) and values on $\mathbb{E}$, and $\mathbb{E}[[\epsilon]]$ stands for the set of formal power series in $\epsilon$ with coefficients in $\mathbb{E}$. We denote by $\mathcal{O}_b(V,\mathbb{E})$ the set of bounded holomorphic functions on $V$  with values in $\mathbb{E}$, for some open set $V$ of $\C^{h}$, for some positive integer $h$. For simplicity, we write $\mathcal{O}(U)$ and $\mathcal{O}_b(V)$ instead of $\mathcal{O}(U,\C)$ and $\mathcal{O}_b(V,\C)$, respectively.

\section{Main problem under study. Analytic solutions}\label{sec1}

In this section, we state the main problem under study, and recall the previous known results associated to the existence and construction of its analytic solutions. All the results and their proofs can be found in detail in~\cite{drelasmal}.

Let $1\le k_1<k_2$ be real numbers, and define $\kappa$ by
$$\frac{1}{\kappa}:=\frac{1}{k_1}-\frac{1}{k_2}.$$
We also fix $q>1$ and consider integers $D,D_1,D_2\ge 3$ and positive integers $d_{D_1},d_{D_2}$. For each $1\le j\le D-1$ we also fix nonnegative integers $\Delta_j$ and positive integers $d_j,\delta_j$. Let us arrange the previous elements as follows:
$$\delta_1=1,\quad \delta_j<\delta_{j+1},\quad 1\le j\le D-2.$$
We additionally assume the following condition holds:
\begin{itemize}
\item[(H1)] For every $j\in\{1,\ldots,D-1\}$ one has
$$\Delta_j\ge d_j,\quad \frac{d_{D_1}-1}{\kappa}+\frac{d_j}{k_2}+1\ge \delta_j,\quad \frac{d_j}{k_1}+1\ge \delta_j,\quad \frac{d_{D_2}-1}{k_2}\ge \delta_j-1,$$
together with
$$k_1(d_{D_2}-1)>k_2d_{D_1}.$$
\end{itemize}

In addition to this, we fix a positive number $\mu$ and polynomials $Q,R_j\in\C[X]$, for $j\in\{1,\ldots,D-1\}$, and $R_{D_1},R_{D_2}\in\C[X]$ such that 
$$\mu>\hbox{deg}(R_{D_j})+1$$ 
for $j\in\{1,2\}$. Moreover, the following statement is assumed:
\begin{itemize}
\item[(H2)] 
$$Q(im)\neq0,\quad R_{D_j}(im)\neq0,$$
for $m\in\R$ and $j\in\{1,2\}$. In addition to this, 
$$\hbox{deg}(Q)\ge\hbox{deg}(R_{D_1})=\hbox{deg}(R_{D_2})\ge\max_{\ell=1,\ldots, D-1}\hbox{deg}(R_{\ell}).$$ 
\end{itemize}

The main equation under study is determined by the previous elements and assumptions.

\begin{align}
Q(\partial_z)\sigma_{q,t}u(t,z,\epsilon)=\sum_{j=1}^{2}(\epsilon t)^{d_{D_j}}\sigma_{q,t}^{\frac{d_{D_j}}{k_j}+1}R_{D_j}(\partial_z)u(t,z,\epsilon)\qquad\qquad\qquad\nonumber\\
\qquad\qquad\qquad+\sum_{\ell=1}^{D-1}\epsilon^{\Delta_{\ell}}t^{d_{\ell}}\sigma_{q,t}^{\delta_{\ell}}(c_{\ell}(t,z,\epsilon)R_{\ell}(\partial_z)u(t,z,\epsilon))+\sigma_{q,t}f(t,z,\epsilon).\label{epral}
\end{align}

Let $\beta>0$. The coefficients $c_{\ell}$, for $\ell\in\{1,\ldots,D-1\}$ and $f$ are defined as an inverse Fourier transform (see Appendix~\ref{secannex}):
$$c_{\ell}(t,z,\epsilon)=\mathcal{F}^{-1}\left(m\mapsto \sum_{p\ge0}C_{\ell,p}(m,\epsilon)(\epsilon t)^p\right)(z),$$
and
$$f(t,z,\epsilon)=\mathcal{F}^{-1}\left(m\mapsto \sum_{p\ge0}F_p(m,\epsilon)(\epsilon t)^p\right)(z)
,$$
respectively. In the previous definitions, $C_{\ell,p}$ turns out to be a holomorphic function with respect to $\epsilon\in D(0,\epsilon_0)$, for some small $0<\epsilon_0<1$, with the function $m\mapsto C_{\ell,p}(m,\epsilon)$ belonging to $\mathcal{C}(\R)$, for all $\epsilon\in D(0,\epsilon_0)$. In addition to this, it satisfies that there exist constants $\Delta_{C,\ell},T_0>0$ such that
$$\sup_{\epsilon\in D(0,\epsilon_0)}|C_{\ell,p}(m,\epsilon)|\le \Delta_{C,\ell}\left(\frac{1}{T_0}\right)^{p}\frac{1}{q^{\frac{p^2\kappa}{2k_1k_2}}} \frac{1}{(1+|m|)^{\mu}}\exp(-\beta|m|),\qquad m\in\R, p\ge0.$$
On the other hand, for all $p\ge0$, the map $F_p$ is a holomorphic function with respect to $\epsilon\in D(0,\epsilon_0)$, and for all $\epsilon\in D(0,\epsilon_0)$ the function $m\mapsto F_{p}(m,\epsilon)$ belongs to $\mathcal{C}(\R)$. Moreover, there exists a constant $\Delta_{F}>0$ such that
$$\sup_{\epsilon\in D(0,\epsilon_0)}|F_{p}(m,\epsilon)|\le \Delta_{F}\left(\frac{1}{T_0}\right)^{p} \frac{1}{(1+|m|)^{\mu}}\exp(-\beta|m|),\qquad m\in\R, p\ge0.$$

From the previous data, it turns out that for every open bounded subset $\mathcal{T}\subseteq\C$ and any fixed $0<\beta'<\beta$, then $c_{\ell}\in\mathcal{O}_b(\mathcal{T}\times H_{\beta'}\times D(0,\epsilon_0))$ for all $\ell\in\{1, \ldots,D-1\}$. 
Here, $H_{\beta}$ stands for the horizontal strip
$$H_{\beta}:=\left\{z\in\C:|\hbox{Im}(z)|<\beta\right\}.$$
In practice, we will fix $\mathcal{T}$ to be a bounded sector with vertex at the origin. In a similar manner, the function $f$ belongs to $\mathcal{O}_b(D(0,\epsilon_0^{-1}T_0/2)\times H_{\beta'}\times D(0,\epsilon_0))$, taking into account the previous considerations.

In~\cite{drelasmal}, it is proved the existence of a family of holomorphic functions of (\ref{epral}) defined on certain domains related to the following construction. We only provide the details involved in the present work, and refer to Definitions 5.1 and 5.2,~\cite{drelasmal}, for further details. For that purpose, we recall the definition of a good covering in $\C^{\star}$.

\begin{defin}\label{defi154}
Let $\varsigma\ge2$ be an integer number. A good covering in $\C^{\star}$ is a family $(\mathcal{E}_p)_{0\le p\le \varsigma-1}$ such that for every $p\in\{0,\ldots,\varsigma-1\}$ the set $\mathcal{E}_p$ is an open sector of positive radius and vertex at the origin. In addition to this, the elements of the good covering satisfy the following statements:
\begin{itemize}
\item $\mathcal{E}_{j}\cap \mathcal{E}_k\neq \emptyset$ for $\{j,k\}\subseteq\{0,\ldots,\varsigma-1\}$ if and only if $|j-k|\le 1$ (by convention, we define $\mathcal{E}_{\varsigma}:=\mathcal{E}_0$). 
\item  There exists a neighborhood of the origin $\mathcal{U}$ such that $\bigcup_{p=0}^{\varsigma-1}\mathcal{E}_p=\mathcal{U}\setminus\{0\}$.
\end{itemize}
\end{defin}

Let us fix a good covering in $\C^{\star}$.

Let $\alpha,\nu\in\R$, with $\nu<1/2$. Associated to the previous good covering, $(\mathcal{E}_p)_{0\le p\le \varsigma-1}$, we consider an open bounded sector $\mathcal{T}\subseteq D(0,r_{\mathcal{T}})$, for some $0<r_{\mathcal{T}}<1$, and with vertex at the origin.
By reducing the radius of the sector, one can always assume that
$$\nu+\frac{k_2}{\log(q)}\log(r_{\mathcal{T}})<0,\quad \alpha+\frac{\kappa}{\log(q)}\log(\epsilon_0r_{\mathcal{T}})<0,\quad \epsilon_0r_{\mathcal{T}}\le q^{\left(\frac{1}{2}-\nu\right)/k_2}/2.$$

We also fix a set of unbounded sectors $(U_{\mathfrak{d}_p})_{0\le p\le \varsigma-1}$, with bisecting direction $\mathfrak{d}_p\in\R$ and a finite family $(\mathcal{R}_{\mathfrak{d}_p}^b)_{0\le p\le \varsigma-1}$ with $\mathcal{R}_{\mathfrak{d}_p}^b=\mathcal{R}_{\mathfrak{d}_p,\tilde{\delta}}\cap D(0,\epsilon_0r_{\mathcal{T}})$, where
\begin{equation}\label{e161}
\mathcal{R}_{\mathfrak{d}_p,\tilde{\delta}}=\left\{T\in\C^{\star}:\left|1+\frac{re^{i\mathfrak{d}_p}}{T}\right|>\tilde{\delta}\hbox{ for all }r\ge0\right\},
\end{equation}
for some $0<\tilde{\delta}<1$. 

For every $p\in\{0,\ldots,\varsigma-1\}$, the direction $\mathfrak{d}_p$ and a positive number $\rho$ are chosen in accordance to the next geometric construction:

\begin{itemize}
\item In view of (H2), it is possible to choose  an unbounded sector 
$$S_{Q,R_{D_1}}=\left\{z\in\C: |z|\ge r_{Q,R_{D_1}},|\hbox{arg}(z)-d_{Q,R_{D_1}}|\le \nu_{Q,R_{D_1}}\right\},$$
for certain $d_{Q,R_{D_1}}\in\R$ and $r_{Q,R_{D_1}},\nu_{Q,R_{D_1}}>0$, such that
$$\frac{Q(im)}{R_{D_1}(im)}\in S_{Q,R_{D_1}},\quad m\in\R.$$
Let $\{q_{0,1}(m),\ldots,q_{d_{D_1}-1,1}(m)\}$ be the roots (in the variable $\tau$) of the polynomial 
$$P_{m,1}(\tau)=\frac{Q(im)}{(q^{1/k_1})^{\frac{k_1(k_1-1)}{2}}}-\frac{R_{D_1}(im)}{(q^{1/k_1})^{\frac{(d_{D_1}+k_1)(d_{D_1}+k_1-1)}{2}}}\tau^{d_{D_1}}.$$

We choose $\mathfrak{d}_p\in\R$ and an infinite sector $U_{\mathfrak{d}_p}$ with bisecting direction $\mathfrak{d}_p\in\R$, among those real directions satisfying
\begin{itemize}
\item[a)] there exists $M_1>0$ such that $|\tau-q_{\ell,1}(m)|\ge M_1(1+|\tau|)$ for all $\ell\in\{0,\ldots,d_{D_1}-1\}$, $m\in\R$ and $\tau\in U_{\mathfrak{d}_p}\cup \overline{D}(0,\rho)$.
\item[b)] there exists $M_2>0$ such that $|\tau-q_{\ell,1}(m)|\ge M_2|q_{\ell,1}(m)|$ for all $\ell\in\{0,\ldots,d_{D_1}-1\}$, $m\in\R$ and $\tau\in U_{\mathfrak{d}_p}\cup \overline{D}(0,\rho)$.
\end{itemize}

For a second geometric condition on the choice of $\mathfrak{d}_p$, we proceed as follows.

\item In view of (H2), it is possible to find an unbounded sector 
$$S_{Q,R_{D_2}}=\left\{z\in\C: |z|\ge r_{Q,R_{D_2}},|\hbox{arg}(z)-d_{Q,R_{D_2}}|\le \nu_{Q,R_{D_2}}\right\},$$
for certain $d_{Q,R_{D_2}}\in\R$ and $r_{Q,R_{D_2}},\nu_{Q,R_{D_2}}>0$, such that
$$\frac{Q(im)}{R_{D_2}(im)}\in S_{Q,R_{D_2}},\quad m\in\R.$$
Let $\{q_{0,2}(m),\ldots,q_{d_{D_2}-1,2}(m)\}$ be the roots (in the variable $\tau$) of the polynomial 
$$P_{m,2}(\tau)=\frac{Q(im)}{(q^{1/k_2})^{\frac{k_2(k_2-1)}{2}}}-\frac{R_{D_2}(im)}{(q^{1/k_2})^{\frac{(d_{D_2}+k_2)(d_{D_2}+k_2-1)}{2}}}\tau^{d_{D_2}}.$$

The direction $\mathfrak{d}_p\in\R$ can also be chosen in such a way that an infinite sector $S_{\mathfrak{d}_p}$ with bisecting direction $\mathfrak{d}_p$ satisfies
\begin{itemize}
\item[a)] there exists $M_{12}>0$ such that $|\tau-q_{\ell,2}(m)|\ge M_{12}(1+|\tau|)$ for all $\ell\in\{0,\ldots,d_{D_2}-1\}$, $m\in\R$ and $\tau\in S_{\mathfrak{d}_p}$.
\item[b)] there exists $M_{22}>0$ such that $|\tau-q_{\ell,2}(m)|\ge M_{22}|q_{\ell,2}(m)|$ for all $\ell\in\{0,\ldots,d_{D_2}-1\}$, $m\in\R$ and $\tau\in S_{\mathfrak{d}_p}$.
\end{itemize}

\item For $p\in\{0,\ldots,\varsigma-1\}$ it holds that $\mathcal{R}_{\mathfrak{d}_p}^{b}\cap \mathcal{R}_{\mathfrak{d}_{p+1}}^{b}\neq\emptyset$ and for every $t\in\mathcal{T}$, $\epsilon\in \mathcal{E}_p$ it holds that $\epsilon t\in \mathcal{R}_{\mathfrak{d}_p}^{b}$ (where $\mathcal{R}_{\mathfrak{d}_\varsigma}^{b}=\mathcal{R}_{\mathfrak{d}_0}^{b}$).
\end{itemize}

The previous geometric construction provides a tuple $\{(\mathcal{R}_{\mathfrak{d}_{p,\tilde{\delta}}})_{0\le p\le \varsigma-1},D(0,\rho),\mathcal{T}\}$ which is known to be associated with the good covering $(\mathcal{E}_{p})_{0\le p\le \varsigma-1}$. In this geometric framework, a family of analytic solutions to (\ref{epral}) can be constructed.


We only give some details on the construction of the analytic solutions which will be needed in the sequel.

\begin{theo}[Theorem 5.3,~\cite{drelasmal}]\label{teo1}
Let $\varsigma\ge2$. In the previous situation regarding the elements involved in the problem (\ref{epral}), and its geometric configuration, assume a good covering in $\C^{\star}$ is given, say $(\mathcal{E}_{p})_{0\le p\le \varsigma-1}$. Let$\{(\mathcal{R}_{\mathfrak{d}_{p,\tilde{\delta}}})_{0\le p\le \varsigma-1},D(0,\rho),\mathcal{T}\}$ be a family associated to the previous good covering. Then, there exist large enough $r_{Q,R_{D_1}},r_{Q,R_{D_2}}>0$, and constants $\varsigma_{c},\varsigma_{F}>0$ such that if 
$$\Delta_{C,\ell}\le \varsigma_{c},\qquad \Delta_{F}\le \varsigma_{F},$$
for all $\ell\in\{1,\ldots,D-1\}$, then for all $p\in\{0,\ldots,\varsigma-1\}$ there exists an analytic solution of (\ref{epral}), $u_p(t,z,\epsilon)\in\mathcal{O}(\mathcal{T}\times H_{\beta'}\times \mathcal{E}_p)$ for every $0<\beta'<\beta$.
\end{theo}

\begin{proof}
Let $p\in\{0,\ldots,\varsigma-1\}$. The construction of the analytic solution is performed in the form of an inverse Fourier and $q-$Laplace transforms by
\begin{equation}\label{esol}
u_p(t,z,\epsilon)=\frac{1}{(2\pi)^{1/2}}\frac{k_2}{\log(q)}\int_{-\infty}^{\infty}\int_{L_{\mathfrak{d}_p}}\frac{w^{\mathfrak{d}_p}_{k_2}(u,m,\epsilon)}{\Theta_{q^{1/k_2}}\left(\frac{u}{\epsilon t}\right)}\frac{du}{u}\exp(izm)dm,
\end{equation}
for the integration path $L_{\mathfrak{d}_p}=[0,\infty)e^{\sqrt{-1}\mathfrak{d}_{p}}$. The function $\Theta_{q^{1/k_2}}(\cdot)$ stands for Jacobi Theta function of order $k_2$ (see Appendix~\ref{secannex}). 

The function $w_{k_2}^{\mathfrak{d}_p}(u,m,\epsilon)$ is obtained from a fixed point argument in adequate Banach spaces of weighted functions. It defines a continuous function on $S_{\mathfrak{d}_p}\times \R\times D(0,\epsilon_0)$ which remains analytic with respect to the first variables on $S_{\mathfrak{d}_p}\times D(0,\epsilon_0)$. Moreover, there exists $\nu\in\R$ and $C_{w_{k_2}^{\mathfrak{d}_p}}>0$ such that
\begin{equation}\label{e218}
\sup_{\epsilon\in D(0,\epsilon_0)}|w_{k_2}^{\mathfrak{d}_p}(u,m,\epsilon)|\le C_{w_{k_2}^{\mathfrak{d}_p}}\frac{1}{(1+|m|)^{\mu}}e^{-\beta|m|}\exp\left(\frac{k_2}{2\log(q)}\log^2|u|+\nu\log|u|\right),
\end{equation}
for all $u\in S_{\mathfrak{d}_{p}}$ and $m\in\R$. The properties of Jacobi Theta function described in Appendix~\ref{secannex} guarantee convergence and holomorphy of (\ref{esol}) in $\mathcal{T}\times H_{\beta'}\times \mathcal{E}_p$, for all $0<\beta'<\beta$.
\end{proof}

In the construction of the analytic solution, the definition of an auxiliary function is needed. We describe its properties in the next result, for the sake of completeness.

\begin{prop}[Proposition 4.2 and Proposition 4.4,~\cite{drelasmal}]\label{prop227}
Assume the hypotheses and constructions involved in Theorem~\ref{teo1} hold. For every $0\le p\le \varsigma-1$ there exists a function $w_{k_1}^{\mathfrak{d}_{p}}$, continuous in $(U_{\mathfrak{d}_{p}}\cup \overline{D}(0,\rho))\times \R\times\mathcal{E}_p$, and holomorphic with respect to its first and third variables on $U_{\mathfrak{d}_{p}}\cup D(0,\rho)$ and $\mathcal{E}_p$ which satisfies the following properties:
\begin{itemize}
\item There exists $C_{p,1}>0$, some $\alpha\in\R$, and some small enough $\delta>0$ such that
$$\sup_{\epsilon\in\mathcal{E}_p}|w_{k_1}^{\mathfrak{d}_{p}}(u,m,\epsilon)|\le C_{p,1}\frac{1}{(1+|m|)^{\mu}}e^{-\beta|m|}\exp\left(\frac{\kappa}{2}\frac{\log^{2}|u+\delta|}{\log(q)}+\alpha\log|u+\delta|\right),$$
for all $u\in U_{\mathfrak{d}_p}\cup\overline{D}(0,\rho)$ and $m\in\R$.
\item The function 
$$\tau\mapsto (\mathcal{L}_{q;1/\kappa}^{\mathfrak{d}_p}(w_{k_1}^{\mathfrak{d}_{p}}(u,m,\epsilon))(\tau)$$
defines a holomorphic function in $\mathcal{R}_{\mathfrak{d}_p,\tilde{\delta}}\cap D(0,r_1)$ for $0<r_1\le q^{\left(\frac{1}{2}-\alpha\right)/\kappa}/2$. In addition to this, it holds that
$$(\mathcal{L}_{q;1/\kappa}^{\mathfrak{d}_p}(w_{k_1}^{\mathfrak{d}_p}(u,m,\epsilon)))(\tau)=w_{k_2}^{\mathfrak{d}_p}(\tau,m,\epsilon),$$
for all $(\tau,m,\epsilon)\in S_{\mathfrak{d}_p}^{b}\times\R\times D(0,\epsilon_0)$, where $S_{\mathfrak{d}_p}^{b}$ is a finite sector of bisecting direction $\mathfrak{d}_p$.
\end{itemize}
\end{prop}

\section{Formal asymptotic results}

In this section, we briefly describe the asymptotic results relating the analytic solutions obtained in Theorem~\ref{teo1} endowed with a different asymptotic representation. First, we recall the asymptotic result attained in~\cite{drelasmal}, associating to such analytic solutions a formal power series expansion in the perturbation parameter. More precisely, the main result obtained in~\cite{drelasmal} is associated to the classical concept of asymptotic expansion in the complex domain, adapted to the framework of $q$-asymptotic expansions, which are more likely to appear when dealing with $q$-difference-equations. Afterwards, we state the main auxiliary results preparing the main result of the work, which is stated in the next section. 

\subsection{First asymptotic approximation}\label{sec31}

The following asymptotic expansion relating the analytic solution and some asymptotic formal representation is considered in this subsection.

\begin{defin}\label{defi1}
Let $\mathcal{E}$ be a bounded sector with vertex at the origin. We fix a complex Banach space $(\mathbb{E},\left\|\cdot\right\|_{\mathbb{E}})$. We also fix a real number $q>1$, and a positive integer $k$. The function $f\in\mathcal{O}(\mathcal{E},\mathbb{E})$ admits the formal power series $\hat{f}(\epsilon)=\sum_{p\ge0}f_p\epsilon^{p}\in\mathbb{E}[[\epsilon]]$ as its $q$-Gevrey asymptotic expansion of order $1/k$ (at the origin) if for every $T\prec\mathcal{E}$ there exist $A,C>0$ such that for all $N\ge0$ one has that
$$\left\|f(\epsilon)-\sum_{p=0}^{N}f_p\epsilon^p\right\|_{\mathbb{E}}\le C A^{N+1}q^{\frac{N(N+1)}{2k}}|\epsilon|^{N+1},$$
for all $\epsilon\in T$.
\end{defin}

The main result in~\cite{drelasmal} reads as follows.

\begin{theo}[Theorem 6.4,~\cite{drelasmal}]\label{teopral0}
Assume the hypotheses of Theorem~\ref{teo1} hold. Let $0<\beta'<\beta$. Then, there exists a formal power series $\hat{u}$ in the form
$$\hat{u}(t,z,\epsilon)=\sum_{p\ge0}h_p(t,z)\frac{\epsilon^{p}}{p!},$$
where $h_p\in\mathcal{O}_{b}(\mathcal{T}\times H_{\beta'})$ satisfies the next features:
\begin{itemize}
\item The formal power series can be written in the form
$$\hat{u}(t,z,\epsilon)=a(t,z,\epsilon)+\hat{u}_1(t,z,\epsilon)+\hat{u}_2(t,z,\epsilon),$$
where $a\in\mathcal{O}_b(\mathcal{T}\times H_{\beta'})\{\epsilon\}$, and $\hat{u}_1,\hat{u}_2\in\mathcal{O}_b(\mathcal{T}\times H_{\beta'})[[\epsilon]]$.
\item For every $p\in\{0,\ldots,\varsigma-1\}$ there exists a splitting of $u$ in the form
$$u(t,z,\epsilon)=a(t,z,\epsilon)+u_{p,1}(t,z,\epsilon)+u_{p,2}(t,z,\epsilon),$$
where $\epsilon\mapsto u_{p,j}(t,z,\epsilon)$ belongs to $\mathcal{O}(\mathcal{E}_p,\mathcal{O}_b(\mathcal{T}\times H_{\beta'}))$, for $j=1,2$ and this function admits $\hat{u}_j(t,z,\epsilon)$ as its $q$-Gevrey asymptotic expansion of order $1/k_j$ on $\mathcal{E}_p$. 
\end{itemize}
In other words, let us write $\hat{u}_j(t,z,\epsilon)=\sum_{\ell\ge0}\tilde{u}_{j,\ell}(t,z)\epsilon^{\ell}\in\mathcal{O}_b(\mathcal{T}\times H_{\beta'})[[\epsilon]]$, for $j=1,2$. Then, for every $p\in\{0,\ldots,\varsigma-1\}$ and given $T_p\prec\mathcal{E}_\ell$, there exist $C,A>0$ such that for all $N\ge0$ one has
$$\sup_{t\in\mathcal{T},z\in H_{\beta'}}\left|u_{p,j}(t,z,\epsilon)-\sum_{\ell=0}^{N}\tilde{u}_{j,\ell}(t,z)\epsilon^{\ell}\right|\le C A^{N+1}q^{\frac{N(N+1)}{2k_j}}|\epsilon|^{N+1},$$
for all $\epsilon\in T_p$.
\end{theo}

\vspace{0.3cm}

\textbf{Remark:} Observe that in the previous result, the Banach space $\mathbb{E}$ of Definition~\ref{defi1} is that of all holomorphic and bounded functions defined in $\mathcal{T}\times H_{\beta'}$, equipped with the norm of the supremum.

\subsection{Preliminary results for a second asymptotic approximation}\label{sec32}

In this subsection, we describe a novel asymptotic representation of the analytic solution to (\ref{epral}). This approach is based on the form of certain families of formal solutions to $q-$difference equations recently discovered in~\cite{carrillolastra}, and the related asymptotic expansion associated to such formal solutions, put forward very recently by H. Tahara in~\cite{tahara24}.

\begin{defin}\label{defi319}
Let $(\mathbb{F},\left\|\cdot\right\|_{\mathbb{F}})$ be a complex Banach space. We fix two bounded sectors with vertex at the origin of the complex plane, $\mathcal{T}$ and $\mathcal{E}$. We also fix a decreasing sequence of positive real numbers $(r_p)_{p\ge0}$ such that
$$\lim_{p\to\infty}r_p=0.$$
Let $f\in\mathcal{O}_b(\mathcal{T}\times\mathcal{E},\mathbb{F})$, and $\hat{f}(t,\epsilon)=\sum_{p\ge0}\varphi_p(t)\epsilon^p$ be a formal power series with $\varphi_p\in\mathcal{O}(\mathcal{T}\cap D(0,r_p),\mathbb{F})$, for every $p\ge0$.

We say that $\hat{f}$ is the 0-Gevrey asymptotic expansion relative to $(r_p)_{p\ge0}$ of $f$ with respect to $\epsilon$ on $\mathcal{E}$ if there exist $C,A>0$ such that for all $N\ge0$ one has
$$\left\|f(t,\epsilon)-\sum_{p=0}^{N}\varphi_p(t)\epsilon^p\right\|_{\mathbb{F}}\le C A^{N+1}|\epsilon|^{N+1},$$ 
for all $t\in\mathcal{T}\cap D(0,r_N)$ and all $\epsilon\in\mathcal{E}$.
\end{defin}

In order to obtain an asymptotic result involving the previous concept, we need to reconsider the main problem under study, and a different manipulation of the analytic solution determined in Theorem~\ref{teo1}. We distinguish two complementary geometric situations in Proposition~\ref{prop293} and Proposition~\ref{prop269}. In these results, we provide information from both, a functional and a sequential point of view. The sequential counterpart is based on the next result.

\begin{lemma}[Lemma 12,~\cite{malek20}]\label{lema293}
Let $k>0$ and $\gamma\in\R$. For every positive integer $N$, the following inequality holds
$$|T|^{-N}|T|^{\gamma}\exp\left(-\frac{k}{2}\frac{\log^2|T|}{\log(q)}\right)\le q^{\frac{\gamma^2}{2k}}(q^{\frac{-\gamma}{k}})^{N}q^{\frac{N^2}{2k}},$$
for $T\in\C^{\star}$.
\end{lemma}

\begin{prop}\label{prop293}
Assume the hypotheses of Theorem~\ref{teo1} hold. Let $p\in\{0,\ldots,\varsigma-1\}$. Assume the infinite sectors $U_{\mathfrak{d}_p}$ and $U_{\mathfrak{d}_{p+1}}$ have nonempty intersection. Then, the following statements hold:
\begin{itemize}
\item[(1)] There exist $K_1>0$ and $K_2\in\R$ such that
\begin{equation}\label{e295}
|u_{p+1}(t,z,\epsilon)-u_{p}(t,z,\epsilon)|\le K_1\exp\left(-\frac{k_2}{2\log(q)}\log^2|\epsilon t|\right)|\epsilon t|^{K_2},
\end{equation}
for every $t\in\mathcal{T}$, $z\in H_{\beta'}$ and $\epsilon\in\mathcal{E}_{p}\cap\mathcal{E}_{p+1}$.
\item[(2)] There exist $K_1>0$ and $K_2\in\R$ such that
\begin{equation}\label{e295b}
|u_{p+1}(t,z,\epsilon)-u_{p}(t,z,\epsilon)|\le K_1q^{\frac{K_2^2}{2k_2}}(q^{-\frac{K_2}{k_2}})^{N}q^{\frac{N^2}{2k_2}}|\epsilon t|^{N},
\end{equation}
for every $t\in\mathcal{T}$, $z\in H_{\beta'}$, $\epsilon\in\mathcal{E}_{p}\cap\mathcal{E}_{p+1}$ and all $N\ge0$.
\end{itemize}
\end{prop}
\begin{proof}
Under the hypotheses made, it is guaranteed that the infinite sector $U_{\mathfrak{d}_{p},\mathfrak{d}_{p+1}}=\{u\in\C^{\star}:\mathfrak{d}_{p}\le\hbox{arg}(u)\le\mathfrak{d}_{p+1}\}$ is contained in $U_{\mathfrak{d}_p}\cap U_{\mathfrak{d}_{p+1}}$. In view of Proposition~\ref{prop227}, the functions $\mathcal{L}_{q;1/\kappa}^{\mathfrak{d}_p}(w_{k_1}^{\mathfrak{d}_p})(\tau,m,\epsilon)$
and $\mathcal{L}_{q;1/\kappa}^{\mathfrak{d}_{p+1}}(w_{k_1}^{\mathfrak{d}_{p+1}})(\tau,m,\epsilon)$ coincide in $(\mathcal{R}^{b}_{\mathfrak{d}_p}\cap \mathcal{R}^{b}_{\mathfrak{d}_{p+1}})\times \R\times D(0,\epsilon_0)$. Therefore, one can define common a function $w_{k_2}^{\mathfrak{d}_p,\mathfrak{d}_{p+1}}(\tau,m,\epsilon)$ on $(\mathcal{R}^{b}_{\mathfrak{d}_{p}}\cup \mathcal{R}^{b}_{\mathfrak{d}_{p+1}})\times \R\times D(0,\epsilon_0)$. This function turns out to be holomorphic with respect to its first variable on $\mathcal{R}^{b}_{\mathfrak{d}_{p}}\cup \mathcal{R}^{b}_{\mathfrak{d}_{p+1}}$, which is the analytic continuation of the previous ones to a common sector. Let us fix $\tilde{\rho}>0$ with $\tilde{\rho}e^{i{\mathfrak{d}}_p}\in\mathcal{R}_{\mathfrak{d}_p}$ and $\tilde{\rho}e^{i{\mathfrak{d}}_{p+1}}\in \mathcal{R}_{\mathfrak{d}_{p+1}}$.
In view of (\ref{esol}), and by the application of Cauchy's theorem, one can deform the integration path in the definition of the analytic solutions (\ref{esol}) to the main problem to arrive at 
\begin{multline*}
u_{p+1}(t,z,\epsilon)-u_{p}(t,z,\epsilon)\\
=\frac{1}{(2\pi)^{1/2}}\frac{k_2}{\log(q)}\int_{-\infty}^{\infty}\int_{L_{\mathfrak{d}_{p+1},\tilde{\rho}}}\frac{w^{\mathfrak{d}_{p+1}}_{k_2}(u,m,\epsilon)}{\Theta_{q^{1/k_2}}\left(\frac{u}{\epsilon t}\right)}\frac{du}{u}\exp(izm)dm\\
-\frac{1}{(2\pi)^{1/2}}\frac{k_2}{\log(q)}\int_{-\infty}^{\infty}\int_{L_{\mathfrak{d}_{p},\tilde{\rho}}}\frac{w^{\mathfrak{d}_{p}}_{k_2}(u,m,\epsilon)}{\Theta_{q^{1/k_2}}\left(\frac{u}{\epsilon t}\right)}\frac{du}{u}\exp(izm)dm\\
+\frac{1}{(2\pi)^{1/2}}\frac{k_2}{\log(q)}\int_{-\infty}^{\infty}\int_{C_{\tilde{\rho},\mathfrak{d}_p,\mathfrak{d}_{p+1}}}\frac{w^{\mathfrak{d}_{p},\mathfrak{d}_{p+1}}_{k_2}(u,m,\epsilon)}{\Theta_{q^{1/k_2}}\left(\frac{u}{\epsilon t}\right)}\frac{du}{u}\exp(izm)dm\\
=I_1-I_2+I_3,
\end{multline*}
with $L_{\mathfrak{d}_j,\tilde{\rho}}=[\tilde{\rho},+\infty)e^{i\mathfrak{d}_j}$ for $j=p,p+1$, and where $C_{\tilde{\rho},\mathfrak{d}_p,\mathfrak{d}_{p+1}}$ is the arc of circle which connects $\tilde{\rho}e^{i\mathfrak{d}_{p}}$ and $\tilde{\rho}e^{i\mathfrak{d}_{p+1}}$ inside $\mathcal{R}^{b}_{\mathfrak{d}_{p}}\cup \mathcal{R}^{b}_{\mathfrak{d}_{p+1}}$. 

We observe from (\ref{e218}) and (\ref{e357}) that
\begin{equation}\label{e313}
|I_1|\le C_{1.1}|\epsilon t|^{1/2}\int_{\tilde{\rho}}^{\infty}\exp\left(\frac{k_2}{2\log(q)}\log^2|u|+\nu\log|u|\right)\exp\left(-\frac{k_2\log^2\left(\frac{|u|}{|\epsilon t|}\right)}{2\log(q)}\right)\frac{d|u|}{|u|^{3/2}},
\end{equation}
for every $(t,z,\epsilon)\in\mathcal{T}\times H_{\beta'}\times(\mathcal{E}_p\cap \mathcal{E}_{p+1})$, and where 
$$C_{1.1}=\frac{C_{w_{k_2}^{\mathfrak{d}_{p+1}}}}{C_{q,k_2}\tilde{\delta}(2\pi)^{1/2}}\frac{k_2}{\log(q)}\int_{-\infty}^{\infty}\frac{1}{(1+|m|)^{\mu}}e^{-\beta|m|-m\hbox{Im}(z)}dm$$
is a positive constant. On the other hand, one has that
\begin{multline*}
\exp\left(\frac{k_2}{2\log(q)}\log^2|u|\right)\exp\left(-\frac{k_2\log^2\left(\frac{|u|}{|\epsilon t|}\right)}{2\log(q)}\right)\\
=\exp\left(\frac{k_2}{2\log(q)}(-\log^2|\epsilon t|+2\log|u|\log|\epsilon t|)\right).
\end{multline*}
We distinguish two situations:
\begin{itemize}
\item[-] If $|u|\ge 1$, then $\log|u|\log|\epsilon t|<0$ due to $|\epsilon|\le \epsilon_0<1$ and $|t|\le r_{\mathcal{T}}<1$, which yields
$$\exp\left(\frac{k_2}{2\log(q)}(2\log|u|\log|\epsilon t|)\right)\le 1.$$
\item[-] If $\tilde{\rho}\le|u|\le 1$, then 
$$\exp\left(\frac{k_2}{2\log(q)}(2\log|u|\log|\epsilon t|)\right)\le |\epsilon t|^{\frac{k_2}{\log(q)}\log(\tilde{\rho})}.$$
\end{itemize}
Therefore, 
$$|I_1|\le C_{1.1}|\epsilon t|^{1/2} \exp\left(-\frac{k_2}{2\log(q)}\log^2|\epsilon t|\right)\left[|\epsilon t|^{\frac{k_2}{\log(q)}\log(\tilde{\rho})}\int_{\tilde{\rho}}^{1}\frac{|u|^{\nu}}{|u|^{3/2}}d|u|+\int_{1}^{\infty}\frac{|u|^{\nu}}{|u|^{3/2}}d|u|\right].$$
Recall that $\nu<1/2$, which leads to the existence of constants $C_{p+1}>0,K_{p+1}\in\R$ such that
\begin{equation}\label{e332}
|I_1|\le C_{p+1}|\epsilon t|^{K_{p+1}} \exp\left(-\frac{k_2}{2\log(q)}\log^2|\epsilon t|\right),
\end{equation}
for every $t\in\mathcal{T}$, $z\in H_{\beta'}$ and $\epsilon\in\mathcal{E}_p\cap\mathcal{E}_{p+1}$.
An analogous reasoning yields the existence of constants $C_{p}>0,K_{p}\in\R$ such that
\begin{equation}\label{e333}
|I_2|\le C_{p}|\epsilon t|^{K_{p}} \exp\left(-\frac{k_2}{2\log(q)}\log^2|\epsilon t|\right),
\end{equation}
for $t\in\mathcal{T}$, $z\in H_{\beta'}$ and $\epsilon\in\mathcal{E}_p\cap\mathcal{E}_{p+1}$.
Finally, regarding (\ref{e357}) we have 
$$|I_3|\le C_{1.1}|\mathfrak{d}_{p+1}-\mathfrak{d}_{p}||\epsilon t|^{1/2}\exp\left(-\frac{k_2\log^2\left(\frac{\tilde{\rho}}{|\epsilon t|}\right)}{2\log(q)}\right)\exp\left(\frac{k_2}{2\log(q)}\log^2(\tilde{\rho})+\nu\log(\tilde{\rho})\right)\frac{1}{\tilde{\rho}^{1/2}},$$
for $t\in\mathcal{T}$, $z\in H_{\beta'}$ and $\epsilon\in\mathcal{E}_p\cap\mathcal{E}_{p+1}$.
Taking into account that 
$$\exp\left(-\frac{k_2\log^2\left(\frac{\tilde{\rho}}{|\epsilon t|}\right)}{2\log(q)}\right)=\exp\left(\frac{-k_2}{2\log(q)}\log^2(\tilde{\rho})\right)|\epsilon t|^{\frac{k_2\log(\tilde{\rho})}{\log(q)}}\exp\left(-\frac{k_2}{2\log(q)}\log^2|\epsilon t|\right),$$
we conclude
\begin{equation}\label{e334}
|I_3|\le C_{p,p+1}|\epsilon t|^{K_{p,p+1}}\exp\left(-\frac{k_2}{2\log(q)}\log^2|\epsilon t|\right),
\end{equation}
with
$$C_{p,p+1}=C_{1.1}|\mathfrak{d}_{p+1}-\mathfrak{d}_{p}|\exp\left(-\frac{k_2}{2\log(q)}\log^2(\tilde{\rho})\right)\exp\left(\frac{k_2}{2\log(q)}\log^2(\tilde{\rho})+\nu\log(\tilde{\rho})\right)\frac{1}{\tilde{\rho}^{1/2}},$$
and
$$K_{p,p+1}=\frac{1}{2}+\frac{k_2\log(\tilde{\rho})}{\log(q)}.$$
From (\ref{e332}), (\ref{e333}) and (\ref{e334}) we conclude (\ref{e295}).

The second part of the proof is a direct consequence of Lemma~\ref{lema293}.
\end{proof}

The complementary situation from that considered in Proposition~\ref{prop293} is studied in Proposition~\ref{prop269}, and is based on the following known result.

\begin{lemma}[Lemma 5.5,~\cite{drelasmal}]\label{lema361}
Assume the hypotheses of Theorem~\ref{teo1} hold. Let $p\in\{0,\ldots,\varsigma-1\}$. If the infinite sectors $U_{\mathfrak{d}_p}$ and $U_{\mathfrak{d}_{p+1}}$ have empty intersection, then there exist $K_{31}>0$ and $K_{41}\in\R$ such that
\begin{multline*}
|\mathcal{L}_{q;1/\kappa}^{\mathfrak{d}_{p+1}}(w^{\mathfrak{d}_{p+1}}_{k_1}(u,m,\epsilon))(\tau)-\mathcal{L}_{q;1/\kappa}^{\mathfrak{d}_{p}}(w^{\mathfrak{d}_{p}}_{k_1}(u,m,\epsilon))(\tau)|\\
\le K_{31}e^{-\beta|m|}\frac{1}{(1+|m|)^{\mu}}\exp\left(-\frac{\kappa}{2\log(q)}\log^2|\tau|\right)|\tau|^{K_{41}},
\end{multline*}
for every $\tau\in(\mathcal{R}_{\mathfrak{d}_p}^b\cap \mathcal{R}_{\mathfrak{d}_{p+1}}^b)$, $m\in\R$ and $\epsilon\in\mathcal{E}_{p}\cap\mathcal{E}_{p+1}$.
\end{lemma}

\begin{prop}\label{prop269}
Assume the hypotheses of Theorem~\ref{teo1} hold. Let $p\in\{0,\ldots,\varsigma-1\}$ and assume the infinite sectors $U_{\mathfrak{d}_p}$ and $U_{\mathfrak{d}_{p+1}}$ have empty intersection. Then, the following statements hold:
\begin{itemize}
\item[(1)] There exist $K_{3}>0$ and $K_{4}\in\R$ such that
\begin{equation}\label{e295c}
|u_{p+1}(t,z,\epsilon)-u_{p}(t,z,\epsilon)|\le K_{3}\exp\left(-\frac{k_1}{2\log(q)}\log^2|\epsilon t|\right)|\epsilon t|^{K_{4}},
\end{equation}
for every $t\in\mathcal{T}$, $z\in H_{\beta'}$ and $\epsilon\in\mathcal{E}_{p}\cap\mathcal{E}_{p+1}$.
\item[(2)] There exist $K_3>0$ and $K_4\in\R$ such that
\begin{equation}\label{e295d}
|u_{p+1}(t,z,\epsilon)-u_{p}(t,z,\epsilon)|\le K_3q^{\frac{K_4^2}{2k_1}}(q^{-\frac{K_4}{k_1}})^{N}q^{\frac{N^2}{2k_1}}|\epsilon t|^{N},
\end{equation}
for every $t\in\mathcal{T}$, $z\in H_{\beta'}$, $\epsilon\in\mathcal{E}_{p}\cap\mathcal{E}_{p+1}$ and all $N\ge0$.
\end{itemize}
\end{prop}
\begin{proof}
We deform the integration path involved in the difference of two consecutive solutions, taking into account the second point in Proposition~\ref{prop227}, in order to arrive at
\begin{multline*}
u_{p+1}(t,z,\epsilon)-u_{p}(t,z,\epsilon)\\
=\frac{1}{(2\pi)^{1/2}}\frac{k_2}{\log(q)}\int_{-\infty}^{\infty}\int_{L_{\mathfrak{d}_{p+1},\tilde{\rho}}}\frac{w^{\mathfrak{d}_{p+1}}_{k_2}(u,m,\epsilon)}{\Theta_{q^{1/k_2}}\left(\frac{u}{\epsilon t}\right)}\frac{du}{u}\exp(izm)dm\\
-\frac{1}{(2\pi)^{1/2}}\frac{k_2}{\log(q)}\int_{-\infty}^{\infty}\int_{L_{\mathfrak{d}_{p},\tilde{\rho}}}\frac{w^{\mathfrak{d}_{p}}_{k_2}(u,m,\epsilon)}{\Theta_{q^{1/k_2}}\left(\frac{u}{\epsilon t}\right)}\frac{du}{u}\exp(izm)dm\\
- \frac{1}{(2\pi)^{1/2}}\frac{k_2}{\log(q)}\int_{-\infty}^{\infty}\int_{C_{\tilde{\rho},\mathfrak{d}_{p+1},1/2(\mathfrak{d}_{p}+\mathfrak{d}_{p+1})}}\frac{w^{\mathfrak{d}_{p}}_{k_2}(u,m,\epsilon)}{\Theta_{q^{1/k_2}}\left(\frac{u}{\epsilon t}\right)}\frac{du}{u}\exp(izm)dm\\
+ \frac{1}{(2\pi)^{1/2}}\frac{k_2}{\log(q)}\int_{-\infty}^{\infty}\int_{C_{\tilde{\rho},\mathfrak{d}_{p},1/2(\mathfrak{d}_{p}+\mathfrak{d}_{p+1})}}\frac{w^{\mathfrak{d}_{p}}_{k_2}(u,m,\epsilon)}{\Theta_{q^{1/k_2}}\left(\frac{u}{\epsilon t}\right)}\frac{du}{u}\exp(izm)dm\\
+\frac{1}{(2\pi)^{1/2}}\frac{k_2}{\log(q)}\int_{-\infty}^{\infty}\int_{L_{0,\tilde{\rho},1/2(\mathfrak{d}_{p}+\mathfrak{d}_{p+1})}}\frac{\mathcal{L}_{q;1/\kappa}^{\mathfrak{d}_{p+1}}(w_{k_1}^{\mathfrak{d}_{p+1}}(u,m,\epsilon))(\tau)-\mathcal{L}_{q;1/\kappa}^{\mathfrak{d}_{p}}(w_{k_1}^{\mathfrak{d}_{p}}(u,m,\epsilon))(\tau)
}{\Theta_{q^{1/k_2}}\left(\frac{\tau}{\epsilon t}\right)}\frac{d\tau}{\tau}\\ 
\hfill\times\exp(izm)dm\\
=I_1-I_2-I_4+I_5+I_6,\\
\end{multline*}
where the integration paths $L_{\mathfrak{d}_{p},\tilde{\rho}}$ and $L_{\mathfrak{d}_{p+1},\tilde{\rho}}$ are as in Proposition~\ref{prop293}, $C_{\tilde{\rho},\mathfrak{d}_{j},1/2(\mathfrak{d}_{p}+\mathfrak{d}_{p+1})}$ is the arc of circle  from $\tilde{\rho}e^{i\mathfrak{d}_{j}}$ to $\tilde{\rho}e^{i\frac{\mathfrak{d}_{p}+\mathfrak{d}_{p+1}}{2}}$ for $j=p,p+1$, and $L_{0,\tilde{\rho},1/2(\mathfrak{d}_{p}+\mathfrak{d}_{p+1})}$ is the segment $[0,\tilde{\rho}]e^{i \frac{\mathfrak{d}_{p}+\mathfrak{d}_{p+1}}{2}}$.

The upper estimates obtained in (\ref{e332}) and (\ref{e333}) are valid now. 

Analogous upper bounds to those determining (\ref{e334}) yield
\begin{equation}\label{e393}
\max\{|I_4|,|I_5|\}\le\frac{\tilde{C}_{p,p+1}}{2}|\epsilon t|^{K_{p,p+1}}\exp\left(-\frac{k_2}{2\log(q)}\log^2|\epsilon t|\right), 
\end{equation}
for all $t\in\mathcal{T}$, $z\in H_{\beta'}$ and $\epsilon\in\mathcal{E}_{p}\cap\mathcal{E}_{p+1}$ and for some constant $\tilde{C}_{p,p+1}>0$.
Finally, Lemma~\ref{lema361} and property (\ref{e357}) of Jacobi Theta function lead us to
\begin{multline*}
|I_6|\le\frac{1}{(2\pi)^{1/2}}\frac{k_2}{\log(q)}K_{31}\frac{1}{C_{q,k_2}\tilde{\delta}}\left(\int_{-\infty}^{\infty}e^{-\beta|m|-m\hbox{Im}(z)}\frac{1}{(1+|m|)^{\mu}}dm\right)|\epsilon t|^{1/2}\\
\times\int_{0}^{\tilde{\rho}}\exp\left(-\frac{\kappa}{2\log(q)}\log^2(r)\right)\exp\left(-\frac{k_2}{2\log(q)}\log^2\left(\frac{r}{|\epsilon t|}\right)\right)r^{K_{41}-\frac{3}{2}}dr.
\end{multline*}
At this point, we observe that 
\begin{multline*}
\exp\left(-\frac{\kappa}{2\log(q)}\log^2(r)\right)\exp\left(-\frac{k_2}{2\log(q)}\log^2\left(\frac{r}{|\epsilon t|}\right)\right)\\
=\exp\left(-\frac{\kappa}{2\log(q)}\log^2(r)\right)\exp\left(-\frac{k_2}{2\log(q)}\log^2(r)-\frac{k_2}{2\log(q)}\log^2|\epsilon t|+\frac{k_2}{\log(q)}\log(r)\log|\epsilon t|\right).$$
\end{multline*}
Therefore,
$$|I_6|\le \tilde{K}_{31}|\epsilon t|^{1/2}\exp\left(-\frac{k_2}{2\log(q)}\log^2|\epsilon t|\right)\int_{0}^{\tilde{\rho}} \exp\left(-\frac{(\kappa+k_2)}{2\log(q)}\log^2(r)\right)r^{\frac{k_2}{\log(q)}\log|\epsilon t|+K_{41}-\frac{3}{2}}dr,$$
with 
$$\tilde{K}_{31}=\frac{1}{(2\pi)^{1/2}}\frac{k_2}{\log(q)}K_{31}\frac{1}{C_{q,k_2}\tilde{\delta}}\left(\int_{-\infty}^{\infty}e^{-\beta|m|-m\hbox{Im}(z)}\frac{1}{(1+|m|)^{\mu}}dm\right).$$

Usual computations guarantee that for every $m_1\in\R$ and $m_2>0$, the function $H(x)=x^{m_1}\exp(-m_2\log^2(x))$ attains its maximum at $x_0=\exp(\frac{1}{2}\frac{m_1}{m_2})$, with $H(x_0)=\exp(\frac{m_1^2}{4m_2})$. Thus,
\begin{multline*}
\int_{0}^{\tilde{\rho}} \exp\left(-\frac{(\kappa+k_2)}{2\log(q)}\log^2(r)\right)r^{\frac{k_2}{\log(q)}\log|\epsilon t|+K_{41}-\frac{3}{2}}dr\le \int_{0}^{\tilde{\rho}}\exp\left(\frac{(\frac{k_2}{\log(q)}\log|\epsilon t|+K_{41}-\frac{3}{2})^2}{4\frac{\kappa+k_2}{2\log(q)}}\right)dr\\
\le \tilde{\rho}K_{32}\exp\left(\frac{k_2^2}{2(\kappa+k_2)\log(q)}\log^2|\epsilon t|\right)|\epsilon t|^{K_{42}},
\end{multline*}
with 
$$K_{42}=\frac{k_2}{\kappa+k_2}\left(K_{41}-\frac{3}{2}\right),\quad K_{32}=\exp\left(\frac{\log(q)}{2(\kappa+k_2)}\left(K_{41}-\frac{3}{2}\right)^2\right).$$
We conclude that 
$$|I_6|\le \tilde{K}_{31}\tilde{\rho}K_{32}|\epsilon t|^{1/2+K_{42}}\exp\left(-\frac{k_2}{2\log(q)}\log^2|\epsilon t|\right)\exp\left(\frac{k_2^2}{2(\kappa+k_2)\log(q)}\log^2|\epsilon t|\right).$$
Taking $K_{3}=\tilde{K}_{31}\tilde{\rho}K_{32}$ and $K_4=1/2+K_{42}$, and bearing in mind that 
$$-k_2+\frac{k_2^2}{\kappa+k_2}=-k_1,$$
we conclude
\begin{equation}\label{e418}
|I_6|\le K_{3}|\epsilon t|^{K_{4}}\exp\left(-\frac{k_1}{2\log(q)}\log^2|\epsilon t|\right).
\end{equation}
The first part of the result follows from (\ref{e332}), (\ref{e333}), (\ref{e393}), (\ref{e418}), and the fact that $k_1<k_2$.

The second part of the proof is a direct consequence of Lemma~\ref{lema293}.
\end{proof}

\section{Second asymptotic approximation of the analytic solutions}\label{secpral}

In this section, we provide a novel sequential version of Ramis-Sibuya theorem in two levels. This is the adaptation to the several level framework of that put forward in Proposition 7.4,~\cite{tahara24}, which is included here for the sake of completeness, and adapted to our needs.

\begin{theo}[Proposition 7.4,~\cite{tahara24}]\label{teorst}
Let $(\mathbb{F},\left\|\cdot\right\|_{\mathbb{F}})$ be a complex Banach space, and $\varsigma\ge 2$ be an integer number. We also fix a bounded sector $\mathcal{T}$ with vertex at the origin. Let $(r_{p})_{p\ge0}$ be a decreasing sequence of positive real numbers, and such that
$$\lim_{p\to\infty}r_{p}=0.$$ 
We also fix a good covering in $\C^{\star}$, say $(\mathcal{E}_p)_{0\le p\le \varsigma-1}$. For every $p\in\{0,\ldots,\varsigma-1\}$, we consider a function $G_p\in\mathcal{O}_b(\mathcal{T}\times \mathcal{E}_p,\mathbb{F})$ (we denote $G_\varsigma:=G_0$). 

Assume that for all $p\in\{0,\ldots,\varsigma-1\}$ there exist $C,H>0$ such that for all $N\ge0$
$$\left\|G_{p+1}(t,\epsilon)-G_{p}(t,\epsilon)\right\|_{\mathbb{F}}\le CH^{N}|\epsilon|^{N}$$
is valid for all $(t,\epsilon)\in(\mathcal{T}\cap D(0,r_{N}))\times (\mathcal{E}_p\cap\mathcal{E}_{p+1})$.

Then, there exists a formal power series
$$\hat{G}(t,\epsilon)=\sum_{p\ge0}\varphi_{p}(t)\epsilon^p,\hbox{ with }\varphi_{p}\in\mathcal{O}(\mathcal{T}\cap D(0,r_{p}),\mathbb{F}),$$
with the next features:
\begin{itemize}
\item for all $0\le p\le \varsigma-1$, there exists bounded holomorphic maps $\Psi_p:\mathcal{T}\times\mathcal{E}_p\to\mathbb{F}$ which admit $\hat{G}(t,\epsilon)$ as $0$-Gevrey asymptotic expansion relative to $(r_p)_{p\ge0}$ with respect to $\epsilon$ and such that 
$$\Delta_{p}(t,\epsilon)=G_{p+1}(t,\epsilon)-G_p(t,\epsilon)=\Psi_{p+1}(t,\epsilon)-\Psi_p(t,\epsilon),$$
for every $(t,\epsilon)\in\mathcal{T}\times (\mathcal{E}_{p}\cap\mathcal{E}_{p+1})$ (with the convention $\Psi_{\varsigma}(t,\epsilon)=\Psi_0(t,\epsilon)$).
\item For all $0\le p\le \varsigma-1$, the map $G_p(t,\epsilon)$ has $\hat{G}(t,\epsilon)$ as $0$-Gevrey asymptotic expansion relative to $(r_p)_{p\ge0}$ with respect to $\epsilon$ on $\mathcal{E}_p$.
\end{itemize}
\end{theo}

The classical version of Ramis-Sibuya theorem can be found in~\cite{hssi}, Lemma XI-2-6. There also exists a version of this classical result involving several levels in the framework of $q-$Gevrey expansions (as mentioned in Definition~\ref{defi1}), which was stated in~\cite{lama15}. The following is a novel sequential version of (RS) Theorem.

\begin{theo}[Multilevel sequential Ramis-Sibuya (RS) Theorem]

Let $(\mathbb{F},\left\|\cdot\right\|_{\mathbb{F}})$ be a complex Banach space. We also fix a bounded sector $\mathcal{T}$ with vertex at the origin. Let $1\le k_1<k_2$ and $\varsigma\ge2$ be an integer number, with $(r_{p,1})_{p\ge0}$ and $(r_{p,2})_{p\ge0}$ be two decreasing sequences of positive real numbers, and
$$\lim_{p\to\infty}r_{p,j}=0,\quad j=1,2.$$

We also fix a good covering in $\C^{\star}$, say $(\mathcal{E}_p)_{0\le p\le \varsigma-1}$. For every $p\in\{0,\ldots,\varsigma-1\}$ we consider $G_p\in\mathcal{O}_b(\mathcal{T}\times \mathcal{E}_p,\mathbb{F})$ and we set a cocycle
$$\Delta_p(t,\epsilon):=G_{p+1}(t,\epsilon)-G_{p}(t,\epsilon),\quad (t,\epsilon)\in\mathcal{T}\times (\mathcal{E}_p\cap\mathcal{E}_{p+1}),$$
with the convention that $G_{\varsigma}:=G_0$. In addition to this, we assume the following properties hold:
\begin{itemize}
\item[(1)] For every $p\in\{0,\ldots,\varsigma-1\}$ the function $G_p$ is bounded on $\mathcal{T}\times\mathcal{E}_p$.
\item[(2)] There exist two nonempty subsets $I_1,I_2\subseteq\{0,\ldots,\varsigma-1\}$, with $I_1\cup I_2=\{0,\ldots,\varsigma-1\}$, and $I_1\cap I_2=\emptyset$ such that:
\begin{itemize}
\item[-] For every $p\in I_1$ there exist $C_1,H_1>0$ such that for all $N\ge0$
$$\left\|\Delta_p(t,\epsilon)\right\|_{\mathbb{F}}\le C_1H_1^{N}|\epsilon|^{N}$$
is valid for all $(t,\epsilon)\in(\mathcal{T}\cap D(0,r_{N,1}))\times (\mathcal{E}_p\cap\mathcal{E}_{p+1})$.
\item[-] For every $p\in I_2$ there exist $C_2,H_2>0$ such that for all $N\ge0$
$$\left\|\Delta_p(t,\epsilon)\right\|_{\mathbb{F}}\le C_2H_2^{N}|\epsilon|^{N}$$
is valid for all $(t,\epsilon)\in(\mathcal{T}\cap D(0,r_{N,2}))\times (\mathcal{E}_p\cap\mathcal{E}_{p+1})$.
\end{itemize}
\end{itemize}
Then, for every $p\in\{0,\ldots,\varsigma-1\}$, there exist $a(t,\epsilon)\in\mathcal{O}_b(\mathcal{T}\times \mathcal{U},\mathbb{F})$, for $\mathcal{U}=\left(\bigcup_{0\le p\le \varsigma-1}\mathcal{E}_p\right)\cup\{0\}$, and for every $p\in\{0,\ldots,\varsigma-1\}$, $G_p^{j}\in\mathcal{O}_b(\mathcal{T}\times\mathcal{E}_p,\mathbb{F})$, for $j=1,2$, and two formal power series
$$\hat{G}^{j}(t,\epsilon)=\sum_{p\ge0}\varphi_{p}^{j}(t)\epsilon^p,\hbox{ with }\varphi_{p}^{j}\in\mathcal{O}(\mathcal{T}\cap D(0,r_{p,j}),\mathbb{F}),$$
for $j=1,2$, such that:
\begin{itemize}
\item $G_p^{j}(t,\epsilon)$ has $\hat{G}^j(t,\epsilon)$ as 0-Gevrey asymptotic expansion relative to $(r_{p,j})_{p\ge0}$ with respect to $\epsilon$ on $\mathcal{E}_p$, for $j=1,2$.
\item For every $p\in\{0,\ldots,\varsigma-1\}$, the function $G_p$ admits the following decomposition
\begin{equation}\label{e503}
G_p(t,\epsilon)=a(t,\epsilon)+G_{p}^{1}(t,\epsilon)+G_{p}^{2}(t,\epsilon).
\end{equation}
\end{itemize}
\end{theo}
\begin{proof}
For every $p\in\{0,\ldots,\varsigma-1\}$, we define the cocycles 
$$\Delta_p^1(t,\epsilon):=\left\{ \begin{array}{ll} \Delta_p(t,\epsilon) & \hbox{ if }p\in I_1 \\ 0 & \hbox{ if }p\in I_2 \end{array} \right.\qquad \Delta_p^2(t,\epsilon):=\left\{ \begin{array}{ll}
0 & \hbox{ if }p\in I_1\\ \Delta_p(t,\epsilon) & \hbox{ if }p\in I_2 \end{array} \right.,
$$
for $(t,\epsilon)\in\mathcal{T}\times (\mathcal{E}_{p}\cap\mathcal{E}_{p+1})$.

As a consequence of Theorem~\ref{teorst}, there exists a bounded holomorphic map $\Psi_p^1:\mathcal{T}\times \mathcal{E}_p\to\mathbb{F}$ such that for all $p\in\{0,\ldots,\varsigma-1\}$
$$\Delta_p^1(t,\epsilon)=\Psi_{p+1}^1(t,\epsilon)-\Psi_{p}^1(t,\epsilon),$$ 
for every $(t,\epsilon)\in\mathcal{T}\times(\mathcal{E}_p\cap\mathcal{E}_{p+1})$. There also exists a family of functions $(\varphi_j^1)_{j\ge0}$, with $\varphi_j^1\in\mathcal{O}_b(\mathcal{T}\cap D(0,r_{j,1}),\mathbb{F})$ for every $j\ge0$. For all $p\in\{0,\ldots,\varsigma-1\}$, the function $\Psi^1_p$ admits the formal power series $\sum_{j\ge0}\varphi_j^1(t)\epsilon^j$ as its 0-Gevrey asymptotic expansion  relative to $(r_{j,1})_{j\ge0}$ with respect to $\epsilon$ on $\mathcal{E}_p$.

In a similar way, there exists a bounded holomorphic map $\Psi_p^2:\mathcal{T}\times \mathcal{E}_p\to\mathbb{F}$ such that for all $p\in\{0,\ldots,\varsigma-1\}$ 
$$\Delta_p^2(t,\epsilon)=\Psi_{p+1}^2(t,\epsilon)-\Psi_{p}^2(t,\epsilon),$$ 
for every $(t,\epsilon)\in\mathcal{T}\times(\mathcal{E}_p\cap\mathcal{E}_{p+1})$. There also exists a family of functions $(\varphi_j^2)_{j\ge0}$, with $\varphi_j^2\in\mathcal{O}_b(\mathcal{T}\cap D(0,r_{j,2}),\mathbb{F})$ for every $j\ge0$. for all $p\in\{0,\ldots,\varsigma-1\}$, the function $\Psi^2_p$ admits the formal power series $\sum_{j\ge0}\varphi_j^2(t)\epsilon^j$ as its 0-Gevrey asymptotic expansion  relative to $(r_{j,2})_{j\ge0}$ with respect to $\epsilon$ on $\mathcal{E}_p$.

We define, for all $p\in\{0,\ldots,\varsigma-1\}$ the function
$$a_p(t,\epsilon):=G_p(t,\epsilon)-\Psi_{p}^1(t,\epsilon)-\Psi_{p}^2(t,\epsilon),$$
which turns out to be a holomorphic and bounded function defined on $\mathcal{T}\times\mathcal{E}_p$, with values in $\mathbb{F}$. We observe that
\begin{multline*}
a_{p+1}(t,\epsilon)-a_{p}(t,\epsilon)=G_{p+1}(t,\epsilon)-\Psi_{p+1}^1(t,\epsilon)-\Psi_{p+1}^2(t,\epsilon)-\left(G_p(t,\epsilon)-\Psi_{p}^1(t,\epsilon)-\Psi_{p}^2(t,\epsilon)\right)\\
=G_{p+1}(t,\epsilon)-G_{p}(t,\epsilon)-\Delta_p^1(t,\epsilon)-\Delta_p^2(t,\epsilon)\\
=G_{p+1}(t,\epsilon)-G_{p}(t,\epsilon)-\Delta_p(t,\epsilon)=0,
\end{multline*}
for all $t\in\mathcal{T}$ and $\epsilon\in\mathcal{E}_p\cap\mathcal{E}_{p+1}$, and all $p\in\{0,\ldots \varsigma-1\}$ (as usual, $a_\varsigma$ is identified with $a_0$). This guarantees the existence of a holomorphic function $a(t,\epsilon)$, defined on $\mathcal{T}\times\left(\cup_{0\le p\le \varsigma-1}\mathcal{E}_p\right)$ and values in $\mathbb{F}$, by $a(t,\epsilon):=a_p(t,\epsilon)$, for all $(t,\epsilon)\in\mathcal{T}\times\mathcal{E}_{p}$. Observe that $\cup_{0\le p\le \varsigma-1}\mathcal{E}_p$ provides a punctured neighborhood of the origin. In addition to this, from the definition of the function $a$ one has that it is a bounded function at the origin, so $a$ indeed defines a holomorphic function on $\mathcal{T}\times \mathcal{U}$, with $\mathcal{U}$ being some neighborhood of the origin.

We have achieved the decomposition (\ref{e503}), with $G_p^j(t,\epsilon):=\Psi_{p}^j(t,\epsilon)$, for $j=1,2$

\end{proof}

The main result of the present work states the second asymptotic representation of the analytic solutions to the main problem (\ref{epral}).

\begin{theo}\label{teopral}
Assume the hypotheses of Theorem~\ref{teo1} hold. Let $0<\beta'<\beta$. Then, there exists a formal power series $\hat{u}(t,z,\epsilon)$ satisfying the next properties:
\begin{itemize}
\item The formal power series $\hat{u}(t,z,\epsilon)=\sum_{p\ge0}h_p(t,z)\frac{\epsilon^p}{p!}$ can be split in the form
$$\hat{u}(t,z,\epsilon)=\mathfrak{a}(t,z,\epsilon)+\hat{\mathfrak{u}}_1(t,z,\epsilon)+ \hat{\mathfrak{u}}_2(t,z,\epsilon),$$
where $\mathfrak{a}(t,z,\epsilon)\in\mathcal{O}_b(\mathcal{T}\times H_{\beta'})\{\epsilon\}$, and for $j=1,2$, the expression
$$\hat{\mathfrak{u}}_j(t,z,\epsilon)=\sum_{p\ge0}h_p^j(t,z)\frac{\epsilon^p}{p!}$$
is a formal power series in $\epsilon$, whose coefficients $t\mapsto h_p^{j}(t,z)$ are holomorphic and bounded functions on $\mathcal{T}\cap D(0,r_p^j)$, with values in the Banach space $\mathcal{O}_b(H_{\beta'})$. The sequence $(r_{p}^{j})_{p\ge0}$ is given by $r_p^j=q^{-\frac{p}{2k_j}}$, for $j=1,2$ and all $p\ge0$.
\item For all $0\le p\le \varsigma-1$, the analytic solution $u_p$ constructed in (\ref{esol}) can also be split in the form
$$u_p(t,z,\epsilon)=\mathfrak{a}(t,z,\epsilon)+\mathfrak{u}_{p,1}(t,z,\epsilon)+\mathfrak{u}_{p,2}(t,z,\epsilon).$$
For $p\in\{0,\ldots,\varsigma-1\}$ and $j=1,2$, the function $(t,\epsilon)\mapsto \mathfrak{u}_{p,j}(t,z,\epsilon)$, as a function with values in the Banach space of holomorphic and bounded functions in $H_{\beta'}$, admits $(t,\epsilon)\mapsto \hat{\mathfrak{u}}_j(t,z,\epsilon)$ as its 0-Gevrey asymptotic expansion relative to $(r_{p}^{j})_{p\ge0}$ on $\mathcal{E}_p$.
\end{itemize}
\end{theo}
\begin{proof}
Let $0<\beta'<\beta$. We also fix $\mathbb{F}$ to be the Banach space of holomorphic and bounded functions on $H_{\beta'}$, endowed with the sup. norm.

We split the set $\{0,\ldots,\varsigma-1\}$ into $I_1$ and $I_2$, with 
$$I_1=\left\{j\in \{0,\ldots,\varsigma-1\}: U_{\mathfrak{d}_j}\cap U_{\mathfrak{d}_{j+1}}\neq\emptyset\right\},$$
and
$$I_2:=\{0,\ldots,\varsigma-1\}\setminus I_1.$$
For every $p\in\{0,\ldots,\varsigma-1\}$, we define the function $(t,\epsilon)\mapsto G_p(t,\epsilon)$ by $G_p(t,\epsilon):H_{\beta'}\to\C$, being
$$(G_p(t,\epsilon))(z)=u_p(t,z,\epsilon),\qquad (t,z,\epsilon)\in\mathcal{T}\times H_{\beta'}\times \mathcal{E}_p.$$
It turns out that $G_p\in\mathcal{O}_b(\mathcal{T}\times\mathcal{E}_p,\mathbb{F})$ for all $p\in\{0,\ldots,\varsigma-1\}$. 

In view of (\ref{e295b}) in Proposition~\ref{prop293}, we have that for all $p\in I_1$ there exist $K_1>0$, $K_2\in\R$ such that
$$\left\|u_{p+1}(t,z,\epsilon)-u_{p}(t,z,\epsilon)\right\|_{\mathbb{F}}\le K_1q^{\frac{K_2^2}{2k_2}}(q^{-\frac{K_2}{k_2}})^{N}|\epsilon|^{N},$$
for every $t\in\mathcal{T}\cap D(0,r_{N}^{2})$, $\epsilon\in\mathcal{E}_p\cap\mathcal{E}_{p+1}$ and all $N\ge0$.

On the other hand, (\ref{e295d}) in Proposition~\ref{prop269} guarantees that for all $p\in I_2$ there exist $K_3>0$ and $K_4\in\R$ such that  
$$\left\|u_{p+1}(t,z,\epsilon)-u_{p}(t,z,\epsilon)\right\|_{\mathbb{F}}\le K_3q^{\frac{K_4^2}{2k_1}}(q^{-\frac{K_4}{k_1}})^{N}|\epsilon|^{N},$$
for every $t\in\mathcal{T}\cap D(0,r_{N}^{1})$, $\epsilon\in\mathcal{E}_p\cap\mathcal{E}_{p+1}$ and all $N\ge0$.

The hypotheses of Multilevel sequential Ramis-Sibuya (RS) Theorem are satisfied with $C_1=K_1q^{\frac{K_2^2}{2k_2}}$, $C_2=K_3q^{\frac{K_4^2}{2k_1}}$, $H_1=q^{-\frac{K_2}{k_2}}$ and $H_2=q^{-\frac{K_4}{k_1}}$. This entails the existence of $\mathfrak{a}_1(t,\epsilon)$ holomorphic and bounded on $\mathcal{T}\times D(0,\epsilon_1)$, for some $\epsilon_1>0$ with values on $\mathbb{F}$; for $j=1,2$ and $p\in\{0,\ldots,\varsigma-1\}$ some function $G_{p}^{j}\in\mathcal{O}_b(\mathcal{T}\times\mathcal{E}_p,\mathbb{F})$ and two formal power series $\hat{G}^j(t,\epsilon)=\sum_{p\ge0}\varphi_p^j(t)\epsilon^p$, with $\varphi_p^j\in\mathcal{O}(\mathcal{T}\cap D(0,r_p^j),\mathbb{F})$ with $G_p^j(t,\epsilon)$ admitting $\hat{G}^j(t,\epsilon)$ as their 0-Gevrey asymptotic expansion relative to $(r_{p}^{j})_{p\ge0}$ with respect to $\epsilon$ on $\mathcal{E}_{p}$, for $j=1,2$. The function $G_p(t,\epsilon)$ is written as $\mathfrak{a}_1(t,\epsilon)+G_p^1(t,\epsilon)+G_p^2(t,\epsilon)$.

The proof is concluded by defining $\mathfrak{a}(t,z,\epsilon):=(\mathfrak{a}_1(t,\epsilon))(z)$, for all $(t,z,\epsilon)\in \mathcal{T}\times H_{\beta'}\times D(0,\epsilon_1)$, and for all $j=1,2$ and $p\in\{0,\ldots,\varsigma-1\}$ we put $\mathfrak{u}_{p,j}(t,z,\epsilon)=(G_p^j(t,\epsilon))(z)$, $h_p^j=\varphi_p^j$, and $\hat{\mathfrak{u}}_j(t,z,\epsilon)=(\hat{G}^j(t,\epsilon))(z)$.
\end{proof}

\begin{corol}
Assume the hypotheses of Theorem~\ref{teopral} hold. Let $0<\beta'<\beta$. Then, for every $p\in\{0,\ldots,\varsigma-1\}$, the analytic solution of (\ref{epral}), $u_p(t,z,\epsilon)$ admits $\hat{u}(t,z,\epsilon)$ as 0-Gevrey asymptotic expansion relative to $(r_p^1)_{p\ge0}$ with respect to $\epsilon$ on $\mathcal{E}_p$, as functions (resp. coefficients of the formal power series) with values in the Banach space of bounded holomorphic functions on $H_{\beta'}$. 
\end{corol}

\begin{proof}
It is a direct consequence of the fact that $r_p^{1}\le r_p^{2}$ for every $p\ge0$ which guarantees holomorphy of the restriction of $h^2_p$ to $(\mathcal{T}\times D(0,r_p^1))\times H_{\beta'}$, and therefore of $h_p$. 
\end{proof}

\section{Appendix}\label{secannex}

\subsection{Inverse Fourier transform}
Let $\beta>0$ and $\mu>1$. Given $f\in\mathcal{C}(\R)$ such that there exists $C>0$ with
$$|f(m)|\le C \frac{1}{(1+|m|)^{\mu}}\exp(-\beta|m|),\qquad m\in\R,$$
one can define 
$$\mathcal{F}^{-1}(f)(x)=\frac{1}{\sqrt{2\pi}}\int_{-\infty}^{\infty}f(m)\exp(ixm)dm,\quad x\in\R.$$
It holds that $\mathcal{F}^{-1}(f)\in\mathcal{O}(H_{\beta}),$ where $H_{\beta}$ stands for the horizontal strip
$$H_{\beta}=\{z\in\C:|\hbox{Im}(z)|<\beta\}.$$

\subsection{Jacobi Theta function}
Let $q>1$. Jacobi Theta function of order $k>0$ is the function defined by
$$\Theta_{q^{1/k}}(z)=\sum_{p\in\Z}\frac{1}{q^{\frac{p(p-1)}{2k}}}z^p,$$
for every $z\in\C^{\star}$. Jacobi Theta function is a holomorphic function in $\C^{\star}$, with an essential singularity at the origin. Moreover, it is the solution on $\C^{\star}$ of the $q-$difference equation
$$\sigma_{q,z}^{\frac{m}{k}}y=q^{\frac{m(m+1)}{2k}}z^{m}y,$$
for all $m\in\Z$.

The asymptotic behavior of Jacobi Theta function near infinity is related to a $q$-exponential growth of order $k$, i.e. for every $\tilde{\delta}>0$ there exists $C_{q,k}>0$ (independent of $\tilde{\delta}$) such that
\begin{equation}\label{e357}
|\Theta_{q^{1/k}}(z)|\ge C_{q,k}\tilde{\delta}\exp\left(\frac{k}{2}\frac{\log^2|z|}{\log(q)}\right)|z|^{1/2},\quad z\in\{z\in\C^{\star}:\inf_{m\in\Z}|1+zq^{\frac{m}{k}}|>\tilde{\delta}\}.
\end{equation}

\subsection{q-Laplace transform}

In this section, $(\mathbb{F},\left\|\cdot\right\|_{\mathbb{F}})$ denotes a complex Banach space.

Let $U_d$ be an infinite sector with vertex at the origin, and bisecting direction $d\in\R$. Let $\rho>0$ and $k>0$. We also fix $f\in\mathcal{O}(D(0,\rho)\cup U_d,\mathbb{F})$, continuous up to its boundary, and assume the existence of $K>0$ and $\alpha\in\R$ such that the following upper estimates hold
$$\left\|f(z)\right\|_{\mathbb{F}}\le K\exp\left(\frac{k\log^2|z|}{2\log(q)}+\alpha\log|z|\right),$$
for all $z\in U_d$ with $|z|\ge \rho$, together with
$$\left\|f(z)\right\|_{\mathbb{F}}\le K,$$
for $z\in \overline{D}(0,\rho)$. The $q$-Laplace transform of order $k$ of $f$ along direction $d$ is defined by
$$\mathcal{L}_{q;1/k}^{d}(f(z))(T)=\frac{k}{\log(q)}\int_{L_{d}}\frac{f(u)}{\Theta_{q^{1/k}}\left(\frac{u}{T}\right)}\frac{du}{u},$$
where $L_{d}$ stands for the path $[0,\infty)\ni r\mapsto re^{\sqrt{-1}d}$.

\begin{lemma}[Lemma 4 and Proposition 6,~\cite{ma17}]
Let $\tilde{\delta}>0$, and consider $f$ as before. Then, $\mathcal{L}_{q;1/k}^{d}(f(z))\in\mathcal{O}(\mathcal{R}_{d,\tilde{\delta}}\cap D(0,r_1),\mathbb{F})$, for all positive $r_1<q^{\left(\frac{1}{2}-\alpha\right)/k}/2$, and $\mathcal{R}_{d,\tilde{\delta}}$ defined by (\ref{e161}). 
\end{lemma}


\textbf{Aknowledgements:} Both authors are partially supported by the project PID2022-139631NB-I00 of Ministerio de Ciencia e Innovaci\'on, Spain. The first author is partially supported by Ministerio de Ciencia e Innovaci\'on-Agencia Estatal de Investigaci\'on MCIN/AEI/10.13039/501100011033 and the European Union ``NextGenerationEU''/ PRTR, under grant TED2021-129813A-I00.


\begin{thebibliography}{99}
\bibitem{ba2} W. Balser, \emph{Formal power series and linear systems of meromorphic ordinary differential equations.} Universitext. Springer-Verlag, New York, 2000. xviii+299 pp.
\bibitem{carrillolastra} S. A. Carrillo, A. Lastra, $q-$Nagumo norms and formal solutions to singularly perturbed $q-$difference equations,	arXiv:2307.15096 [math.GM] 
\bibitem{dre} T. Dreyfus, \emph{Building meromorphic solutions of q-difference equations using a Borel–Laplace summation}, Int. Math. Res. Not. 15, (2015) 6562--6587.
\bibitem{drelasmal} T. Dreyfus, A. Lastra, S. Malek, \emph{Multiple-scale analysis for some linear partial $q-$difference and differential equations with holomorphic coefficients}, Advances in Difference Equations, 2019:326, 2019.
\bibitem{hssi} P. Hsieh, Y. Sibuya, \emph{Basic theory of ordinary differential equations}. Universitext. Springer-Verlag, New York, 1999.
\bibitem{icmi} K. Ichinobe, S. Michalik, \emph{On the summability and convergence of formal solutions of linear q-difference-differential equations with constant coefficients}. Math. Ann. (2023). 
\bibitem{lama15} A. Lastra, S. Malek, \emph{On parametric multilevel $q-$Gevrey asymptotics for some linear $q-$difference-differential equations}, Advances in Difference Equations (2015) 2015:344.
\bibitem{lamasa0} A. Lastra, S. Malek, J. Sanz, \emph{On $q-$asymptotics for linear $q-$difference-differential equations with Fuchsian and irregular singularities}. J. Differ. Equations 252 (2012), no. 10, 5185--5216.
\bibitem{lami} A. Lastra, S. Michalik, \emph{On sequences preserving q-Gevrey asymptotic expansions}. Anal.Math.Phys. 14, 17 (2024).
\bibitem{loday} M. Loday-Richaud, Divergent series, summability and resurgence. II. Simple and multiple summability. Lecture Notes in Mathematics, 2154. Springer, 2016. 
\bibitem{ma11} S. Malek, \emph{Parametric Gevrey asymptotics for a q-analog of some linear initial value problem}, Funkc. Ekvacioj 60(1), (2017) 21--63.
\bibitem{ma17} S. Malek, \emph{Parametric Gevrey asymptotics for a $q-$analog of some linear initial value problem}, Funkc. Ekvacioj 60 no. 1 (2017), 21--63.
\bibitem{malek20} S. Malek, \emph{On a partial $q-$analog of a singularly perturbed problem with fuchsian and irregular time singularities}, Abstract and Applied Analysis, vol. 2020 (2020) Article ID 7985298. 
\bibitem{prrasp} D. W. Pravica, N. Randriampiry, M. J. Spurr, \emph{On $q-$advanced spherical Bessel functions of the first kind and perturbations of the Haar wavelet.} Appl. Comput. Harmon. Anal. 44, No. 2 (2018), 350-413. 
\bibitem{prrasp1} D. W. Pravica, N. Randriampiry, M. J. Spurr, \emph{Solutions of a class of multiplicatively advanced differential equations.} C. R., Math., Acad. Sci. Paris 356, No. 7 (2018), 776--817. 
\bibitem{tahara24} H. Tahara, \emph{Asymptotic existence theorem for formal power series solutions of singularly perturbed linear q-difference equations}, Math. Ann. (2024). 
\bibitem{tahara2} H. Tahara, \emph{On the summability of formal solutions of some linear $q-$difference-differential equations}, Funkc. Ekvacioj, Ser. Int. 63, No. 2 (2020), 259--291. 
\bibitem{taya2} H. Tahara, H. Yamazawa, \emph{$q-$analogue of summability of formal solutions of some linear $q-$difference-differential equations}, Opuscula Math. 35 (2015), no. 5, 713--738.
\bibitem{vizh} L. Di Vizio, C. Zhang, \emph{On q-summation and confluence}, Ann. Inst. Fourier (Grenoble) 59(1), (2009) 347--392.
\bibitem{ya1} H. Yamazawa, \emph{Holomorphic and singular solutions of $q-$difference-differential equations of Briot-Bouquet type}, Funkcial. Ekvac. 59 (2016), no. 2, 185--197.
\end{thebibliography}
\end{document}